\documentclass{article}

\usepackage{arxiv}

\usepackage{bm}
\usepackage{subcaption}
\usepackage{amsmath}
\usepackage{amsfonts}
\usepackage{listings}

\usepackage[svgnames]{xcolor}

\lstdefinelanguage{XML}
{
  basicstyle=\ttfamily\footnotesize,
  morestring=[b]",
  moredelim=[s][\bfseries\color{Maroon}]{<}{\ },
  moredelim=[s][\bfseries\color{Maroon}]{</}{>},
  moredelim=[l][\bfseries\color{Maroon}]{/>},
  moredelim=[l][\bfseries\color{Maroon}]{>},
  morecomment=[s]{<?}{?>},
  morecomment=[s]{<!--}{-->},
  commentstyle=\color{DarkOliveGreen},
  stringstyle=\color{blue},
  identifierstyle=\color{red}
}

\usepackage{color}

\usepackage{pgfplots}
\usepackage{pgfplotstable}
\pgfplotsset{
compat=newest, 
tick label style={font=\footnotesize}, 
}

\usepackage{overpic}

\usepackage{tikz}
\usetikzlibrary{mindmap,shadows}
\usetikzlibrary{arrows,intersections}
\usepackage{smartdiagram}

\usepackage{booktabs}
\usepackage[]{units}

\title{openCFS-Data: Data Pre-Post-Processing Tool for openCFS}

\author{ Stefan Schoder \\
	Group of Aeroacoustics and Vibroacoustics, IGTE\\
	TU Graz\\
	Inffeldgasse 18, 8010 Graz \\
	\texttt{stefan.schoder@tugraz.at} \\
	\And
	Klaus Roppert \\
	Group of Multiphysics, IGTE\\
	TU Graz\\
	Inffeldgasse 18, 8010 Graz \\
	\texttt{klaus.roppert@tugraz.at} \\
}

\renewcommand{\shorttitle}{openCFS Software}

\usepackage{hyperref}
\hypersetup{
pdftitle={openCFS Software},
pdfsubject={preprint},
pdfauthor={Schoder, Roppert},
pdfkeywords={FEM Sotware, open Source},
}

\begin{document}
\maketitle

\begin{abstract}
	Many numerical simulation tools have been developed and are on the market, but there is still a strong need for appropriate tools capable to simulate multi-field problems, especially in aeroacoustics. Therefore, openCFS provides an open-source framework for implementing partial differential equations using the finite element method. Since 2000, the software has been developed continuously. The result of is openCFS (before 2020 known as CFS++ Coupled Field Simulations written in C++). In this paper, we present for the first time the CFS-Data, the open-source pre-post-processing part of openCFS with a focus on the aeroacoustic source computation (called filters).
\end{abstract}


%

\section{Introduction}
\label{sec:Intro} 
Within this contribution, we concentrate on the openCFS \cite{CFS} module \textit{openCFS-Data}. An alternative data processing tool is the recently developed pyCFS-data \cite{wurzinger2024pycfs}. 

When establishing an XML file for CFS-Data, it is fundamental that the pipeline, existing of different CFS-Data filters, is closed. The pipeline has to start with the step value definition and has to be followed by the input filter and end with the output filter. In between, multiple filters can be added, serial or parallel.

\textbf{Defining Step Value Definition:}
It is possible to define input data for the time and frequency domain. However, not all filters are capable of processing data in the frequency domain. 
\begin{lstlisting}[language=XML]
  <pipeline>
  
     <stepValueDefinition>
      <startStop>
        <startStep value="..."/>
        <numSteps value="..."/> 
        <startTime value="..."/>
        <delta  value="..."/>
        <deleteOffset  value="no"/>
      </startStop>
    </stepValueDefinition>

  </pipeline>   
\end{lstlisting}  

\begin{figure}
    \centering
    \includegraphics[scale=3]{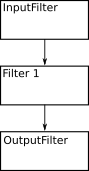}
    \caption{Serial alignment of a filter.}
    \label{fig:serial}
\end{figure}

\begin{figure}
    \centering
    \includegraphics[scale=3]{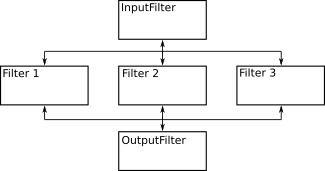}
    \caption{Parallel alignment of filters in a pipeline.}
    \label{fig:multiple}
\end{figure}

Filters can be designed and aligned in a serial or parallel way (see Fig. \ref{fig:serial} and \ref{fig:multiple}).

Such a serial alignment results in the following structure.
\begin{lstlisting}[language=XML]
  <pipeline>
  
     <stepValueDefinition>
      <startStop>
        <startStep value="..."/>
        <numSteps value="..."/> 
        <startTime value="..."/>
        <delta  value="..."/>
        <deleteOffset  value="no"/>
      </startStop>
    </stepValueDefinition>
    
    <meshInput id="inputFilter" gridType="fullGrid"  >      

    </meshInput>
    
    <interpolation type="FieldInterpolation_Cell2Node" id="interp1" inputFilterIds="inputFilter">
    </interpolation>       
    
    <meshOutput id="Outout" inputFilterIds="interp1">

    </meshOutput>
  </pipeline> 
\end{lstlisting}  
whereas a parallel alignment is set up in the following.
\begin{lstlisting}[language=XML]
  <pipeline>
  
     <stepValueDefinition>
      <startStop>
        <startStep value="..."/>
        <numSteps value="..."/> 
        <startTime value="..."/>
        <delta  value="..."/>
        <deleteOffset  value="no"/>
      </startStop>
    </stepValueDefinition>
    
    <meshInput id="inputFilter" gridType="fullGrid"  >      

    </meshInput>

    <interpolation type="FieldInterpolation_Cell2Node" id="interp1" inputFilterIds="inputFilter">
    </interpolation> 
       
    
    <interpolation type="FieldInterpolation_Cell2Node" id="interp2" inputFilterIds="inputFilter">
    </interpolation>    
    
    
    <interpolation type="FieldInterpolation_Cell2Node" id="interp3" inputFilterIds="inputFilter">
    </interpolation>       
    
    <meshOutput id="Outout" inputFilterIds="interp1 interp2 interp3">

    </meshOutput>
  </pipeline>
\end{lstlisting}  
Within the pipeline, different filters can be arranged. The following filter classes are available today:
\begin{itemize}
    \item Interpolation Filters
    \item Conservative Interpolation Filters
    \item Aeroacoustic Source Terms
    \item Synthetic Sources
    \item Data Processing
\end{itemize}
In this contribution, we discuss the input format, output format and field interpolation possibilities (filters) of CFS-Data in more detail.

\section{IO Formats and Definitions}
The data processing tool of openCFS offers the option of importing Ensight-files and hdf5-files (hierarchical data format), whereas the export of mesh-based field data is by default based on hdf5, which is the native data format of openCFS. Additionally, reading of meshes (e.g. target mesh for interpolation) in cgns or cdb format is supported.  Field data can thereby be defined on the nodes or the cell centroids of a computational grid in the time or the frequency domain.

\subsection{Input definition}

The first block of the XML-scheme defines the time domain of the input data to be read. The following XML-snippet illustrates a typical setting.

\begin{lstlisting}[language=XML]
 <stepValueDefinition>
	<startStop>
		<startStep value="0"/>
		<numSteps value="10"/>
		<startTime value="1e-05"/>
		<delta value="1e-05"/>
		<deleteOffset value="no"/>
	</startStop>
</stepValueDefinition>
\end{lstlisting} 
\begin{itemize}
    \item delta: time step size in seconds for data import ("CFS time step". Thereby, the time step size can be a whole multiple of the time step provided by the input file to read every 2nd, 3rd etc. time step. If for example the input data time step is 1s and the defined "CFS time step" is 2s, every second time step is read.
\item startStep: time offset in multiples of the "CFS time step".
\item numSteps: number of time steps to be read.
\item startTime: offset in seconds according to the input data time values.
\item deleteOffset: delete the offset resulting from startTime. If set to yes and startStep to zero, the first time value of the output data will be the "CFS time step" size.
\end{itemize}

The first time step which is read from the input file is the time step corresponding to the time *startStep*  *delta* + *startTime* in seconds. Thereby, the offset resulting from *startTime* can be deleted for the output by enabling the *deleteOffset*-tag.
If data is processed in frequency domain, the same tags are used (*startTime* defines the start frequency in this case).

Subsequently, the mesh-based input data is provided by

\begin{lstlisting}[language=XML]
<meshInput id="input">
	<inputFile>
		<hdf5 fileName="pathToInputFile/InpuFile.hdf5"/>
	</inputFile>
</meshInput>
\end{lstlisting}            
in case of using the hfd5 format (e.g. *openCFS* simulation file).

*Ensight* data is considered in the XML-scheme by

\begin{lstlisting}[language=XML]
<meshInput id="input" gridType="fullGrid">
	<inputFile>
		<ensight fileName="pathToInputFile/InputFile.case" fixFVPyramids="yes" readFVMesh="no">
			<variableList>
				<variable CFSVarName="cfsQuantity1" EnsightVarName="EnsightQuantity1"/>
				<variable CFSVarName="cfsQuantity2" EnsightVarName="EnsightQuantity2"/>
			</variableList>
		</ensight>
	</inputFile>
</meshInput>
\end{lstlisting}         
where in *fileName* the location of the *Ensight* master file (.case or .encas) needs to be provided and the therein defined quantities (e.g. velocity, pressure) need to be defined by *EnsightVarName*. In contrast, the quantities of *hdf5* files are identified automatically by *openCFS*.

\subsection{Output definition}

The processed field data (*resultQuantity1* and *resultQuantitiy2*) in the following XML-snippet) is exported in the native hdf5-format. Thereby, the filename is defined in the XML scheme and the default file extension CFS can be adapted if required. Furthermore, the default compression level of the hdf5-file of 1 can be modified. Furthermore, external files can be enabled, where the field data of each time/frequency step is written to a separate HDF file and the master file includes the mesh data, further file information, and the links to the external files of each time step. To explore hdf5 files and get an understanding of the structure, *HDFview* is recommended. The results of multiple filters (e.g., *filterID1*, *filterID2* in the following example XML) can be either written to all regions or to specified regions of the output mesh as indicated in the snippet.
\begin{lstlisting}[language=XML]
<meshOutput id="OutputFileName" inputFilterIds="filterID1,filterID2">
	<outputFile>
		<hdf5 extension="cfs" compressionLevel="1" externalFiles="no"/>
	</outputFile>
	<saveResults>
		<result resultName="resultQuantity1">
			<allRegions/>
		</result>
		<result resultName="resultQuantity2">
			<regionList>
				<region name="region1"/>
				<region name="region2"/>
			</regionList>
		</result>
	</saveResults>
</meshOutput>
\end{lstlisting} 

It is important if the exported data will be the input of a subsequent *openCFS* simulation, *openCFS* variable names must be used for the declaration of field quantities. Thus, for the acoustic PDE, one of the following names must be chosen.

General acoustic and fluid mechanic quantities:
\begin{itemize}
\item acouPressure
\item acouVelocity
\item acouPotential
\item acoutIntensity
\item fluidMechVelocity
\item meanFluidMechVelocity
\item fluidMechPressure
\item fluidMechDensity
\item fluidMechVorticity
\item fluidMechGradPressure
\end{itemize}

Aeroacoustic Source Terms:
\begin{itemize}
\item acouRhsLoad (general)
\item acouRhsLoadP (PCWE)
\item vortexRhsLoad (Vortex Sound Theorie)
\item acouDivLighthillTensor (Lighthill's acoustic analogy)
\end{itemize}

\section{Interpolation filters} \label{Sec:W}

\subsection{Node2Cell}
The node to cell interpolation filter takes nodal loads and connects them to the cell center, of the cell defined by those nodes
\begin{equation}
e_{\square} = \sum_{i=1}^\mathrm{n}  v_i \, .
\end{equation}
Thereby, $e_\square$ is the load located to the cell, $\mathrm{n}$ the number of nodes of one element, and $v_i$ the nodal loads. 
The following example shows this methodology by considering one tetrahedral element:
\begin{figure}[ht!]
    \centering
    \includegraphics[scale=1.8]{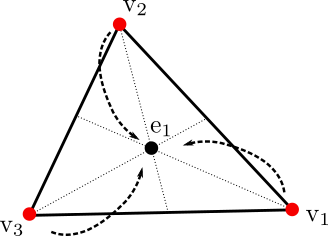}
    \caption{Cell to node interpolator.}
    \label{fig:CC}
\end{figure}
\begin{lstlisting}[language=XML]
    <interpolation type="FieldInterpolation_Cell2Node" id="..." inputFilterIds="...">
      <targetMesh>
        <hdf5 fileName="..."/>
      </targetMesh>
      <singleResult>
        <inputQuantity resultName="..."/>
        <outputQuantity resultName="..."/>
      </singleResult>
      <regions>
        <sourceRegions>
          <region name="..." />
        </sourceRegions>
        <targetRegions>
          <region name="..."/>
        </targetRegions>
      </regions>
    </interpolation>
\end{lstlisting} 
Note it is important that the target mesh is ** not ** the same mesh as the source mesh. Transform source mesh into new, empty mesh by using cfs -g. It is possible to interpolate from a volume onto (curved) surfaces. However, the input data can not be from curved surfaces!

\subsection{Cell2Node}
The cell to node interpolation filter takes element loads and divides it onto the nodes that build the cell.  
\begin{equation}
v_{\square} = \frac{1}{\mathrm{n}} e_i \, .
\end{equation}
The following example shows this methodology by considering one tetrahedral element:
\begin{figure}[ht!]
    \centering
    \includegraphics[scale=1.8]{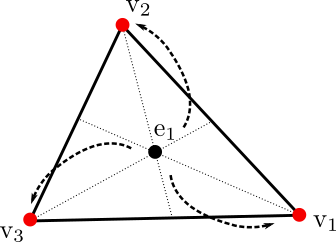}
    \caption{Node to cell interpolator.}
    \label{fig:CC}
\end{figure}
\begin{lstlisting}[language=XML]
    <interpolation type="FieldInterpolation_Node2Cell" id="..." inputFilterIds="...">
      <targetMesh>
        <hdf5 fileName="..."/>
      </targetMesh>
      <singleResult>
        <inputQuantity resultName="..."/>
        <outputQuantity resultName="..."/>
      </singleResult>
      <regions>
        <sourceRegions>
          <region name="..." />
        </sourceRegions>
        <targetRegions>
          <region name="..."/>
        </targetRegions>
      </regions>
    </interpolation>
\end{lstlisting} 
Note it is important that the target mesh is ** not ** the same mesh as the source mesh. Transform source mesh into new, empty mesh by using cfs -g. It is possible to interpolate from a volume onto (curved) surfaces. However, the input data can not be from curved surfaces!

\subsection{Nearest Neighbour}

The nearest neighbor interpolation of openCFS rests upon the inverse distance weighting (Shepard's method). 
\begin{figure}[ht!]
    \centering
    \includegraphics[scale=1.8]{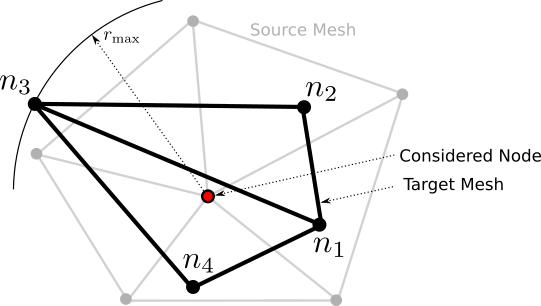}
    \caption{Nearest neighbor interpolator.}
    \label{fig:CC}
\end{figure}
Based on the defined number of neighbours $n$ from the source mesh, the nearest neighbours are searched and their distance $r$ to the considered node is computed. 
Based on this distances the weights $w$ are computed.
\begin{equation}
w_i = \left( \frac{R_{\mathrm{ max}}-r_i}{R_{\mathrm{ max}} r_i} \right)^p \, ,
\end{equation}
with $R_\mathrm{ max} = 1.01 r_{\mathrm{ max}}$ as being 1.01 times of the maximal distance $r_{\mathrm{ max}}$, and $p$ as the interpolation exponent. Shepard stated $1 \leq p \leq 3$. Increasing $ p $ means that values that are further away are taken into account more. Finally, each value of each node $v_1$ that is taken into account is weighted to compute the new value $v_{\mathrm n}$ of the considered node
\begin{equation}
v_{\mathrm n} = \sum_{i=1}^{n} \frac{w_i v_i }{\sum_{i=1}^{n}w_i} \, .
\end{equation}

In CFSDat the Nearest Neighbor interpolation is defined in the xml-file as following:
\begin{lstlisting}[language=XML]
<interpolation type="FieldInterpolation_NearestNeighbour" inputFilterIds="..." id="...">
    <IntSchemeNN interpolationExponent="..." numNeighbours="..." globalFactor="..."/>                   
    <targetMesh>
        <hdf5 fileName="..."/>
    </targetMesh>
    <singleResult>
        <inputQuantity resultName="..."/>
        <outputQuantity resultName="..."/>
    </singleResult>
    <regions>
        <sourceRegions>
            <region name="..."/>
        </sourceRegions>
        <targetRegions>
            <region name="..."/>
        </targetRegions>
    </regions>
</interpolation>
\end{lstlisting} 
*inputFilterIds* must refer to the id of an input filter which contains the defined "sourceRegion" and "quantity". As target mesh, any compatible file input can be used (see chapter Data Input/Output).
\begin{itemize}
    \item interpolationExponent: Exponent for calculation of interpolation weight function.
    \item globalFactor: Global Factor for scaling the result (usually 1).
    \item numNeighbours: Number of considered Nodes.
\end{itemize}

\subsection{Radial Basis Function Interpolation}

This interpolation filter interpolates field data from a source mesh to a target mesh based on radial basis functions (RBFs), where the local Wendland kernel together with a modified Shepard's method was chosen. Thereby, node-based data can be interpolated either to the cell centroids or the nodes of the source mesh, whereas element-based data can currently only be interpolated to the cell-centroids. The theoretical background of the implemented interpolation scheme is published in \cite{schoder2020radial}.

Like for all interpolation filters, the input and output quantitiy, the source and target regions, and the target mesh have to be defined within the xml-scheme. By using the optional tag *noSlipWall*, the output quantity of the nodes on the specified wall region (boundary surface) are set to zero. This tag is intended to consider no-slip-walls for the interpolation of the flow velocity but can theoretically also be used for other quantities. 

The field data is interpolated to the nodes of the target mesh  by default. By setting *useElemAsTarget* *true*, the field data is interpolated to the cell centroids of the target mesh. The tag *interpolationExponent* handles the locality of the approximation on the target mesh and is considered within the weight function (for details see [@schoder2020c]). The larger it is, the more local the approach making it less accurate but capable of resolving stronger gradients. Furthermore, *globalFactor* allows a scaling of the output by the defined factor.

The number of neighbour points $N_q$ (*numNeighbours*) and influence points $N_w$ (*numNeighbours\_weight*) are by default set to $N_q=18$ and $N_w=13$, which is a compromise of accuracy and numerical efficiency for 2D and 3D interpolations. The interopolation result can be enhanced by tuning these parameters specifically for the involved source and target mesh. By setting the tag *useCGAL4RBF* to *true*, the nearest neighbor search will be carried out be by \href{https://www.cgal.org}{CGAL}, which requires the *CGAL 3.3* build (cmake flag must be set in the cmake build properties(../Installation/linuxinstall/)). 

\begin{lstlisting}[language=XML]
<interpolation type="FieldInterpolation_RBF" id="interpolationRBCF" inputFilterIds="input">
	<IntSchemeRBF numNeighbours="..." numNeighbours_weight="..." globalFactor="..." useCGAL4RBF="false"/>
      	<useElemAsTarget>false</useElemAsTarget>
	<noSlipWall name="noSlipWall"/>
	<targetMesh>
		<hdf5 fileName="pathToMeshFile/meshFile.cfs"/>
	</targetMesh>
	<singleResult>
                <inputQuantity resultName="inputQuantity"/>
                <outputQuantity resultName="outputQuantity"/>
	</singleResult>
	<regions>
		<sourceRegions>
			<region name="sourceRegionR" />
		</sourceRegions>
		<targetRegions>
                    <region name="targetRegion"/>
                </targetRegions>
	</regions>
</interpolation>
\end{lstlisting} 

\subsection{Finite Element (FE) based Interpolation}
This filter (beta version) is currently under development and not all functionality is fully tested. A preliminary documentation can be found online in the CFS repository.

\section{Conservative source integration}

The conservative interpolation computes the Right-hand Side force $f$ of a partial differential equation
\begin{equation}
\int_\Omega f \varphi {\mathrm {d}} \Omega \, ,
\end{equation}
using first order nodal FEM to solve the system of equations. If you would like to use the conservative filters for higher order FEM, please contact the authors. \cite{schoder2021application} covers the theory of the implemented conservative interpolation filters. For further development and future publications on the topic, we are happy to collaborate with you.

\paragraph{Pipeline setup}
During the simulation setup, we will discuss how the variables inside the publication relate to the XML-scheme of the simulation setup. Having this knowledge, it should be able for you to start your simulation using conservative interpolators. CGNS or ENSIGHT input data is supported.

\begin{lstlisting}[language=XML]
<?xml version="1.0" encoding="UTF-8"?>
<cfsdat ...>
  <pipeline>

    <stepValueDefinition>
      ...
    </stepValueDefinition>
    
    <meshInput id="inputFilter" gridType="fullGrid">
      ...
    </meshInput>
...
    <interpolation type="TYPE" id="interp" inputFilterIds="inputFilter">
      <targetMesh>
        <hdf5 fileName="../mesh.h5"/>
      </targetMesh>
      <singleResult>
        <inputQuantity resultName="..."/>
        <outputQuantity resultName="..."/>
      </singleResult>
      <regions>
        <sourceRegions>
          <region name="..."/>
        </sourceRegions>
        <targetRegions>
          <region name="..."/>
        </targetRegions>
      </regions>
    </interpolation>
...
    <meshOutput id="acousticSources" inputFilterIds="interp">
      ...
    </meshOutput>

  </pipeline>
</cfsdat>
\end{lstlisting}

There are two variants implemented and their application limits are discussed in \cite{schoder2021application}. The cell centroid interpolator is defined by the type variable 
\begin{lstlisting}[language=XML]
type="FieldInterpolation_Conservative_CellCentroid"
\end{lstlisting} 
based on figure \ref{fig:CC}, the interpolation is carried out.
\begin{figure}[ht!]
    \centering
    \includegraphics[scale=0.3]{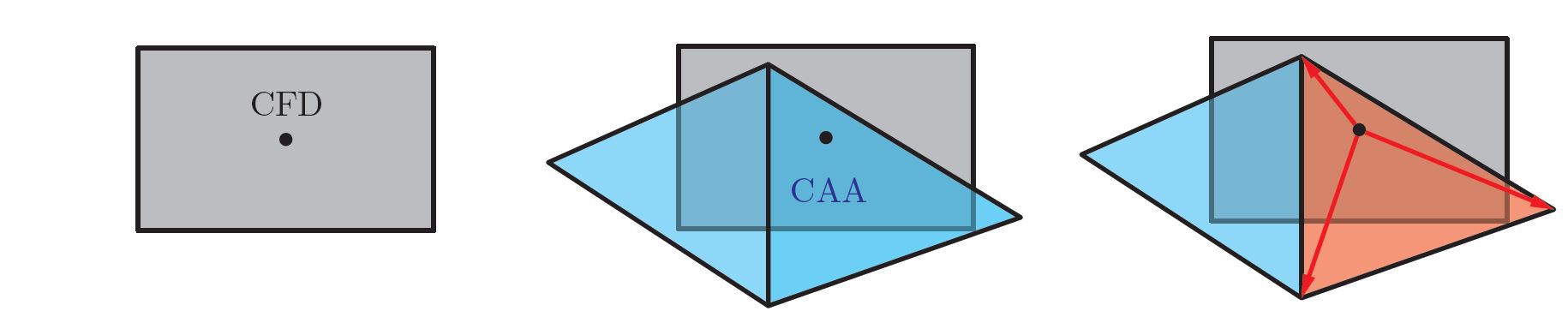}
    \caption{Cell centroid interpolator.}
    \label{fig:CC}
\end{figure}

\begin{figure}[ht!]
    \centering
    \includegraphics[scale=0.3]{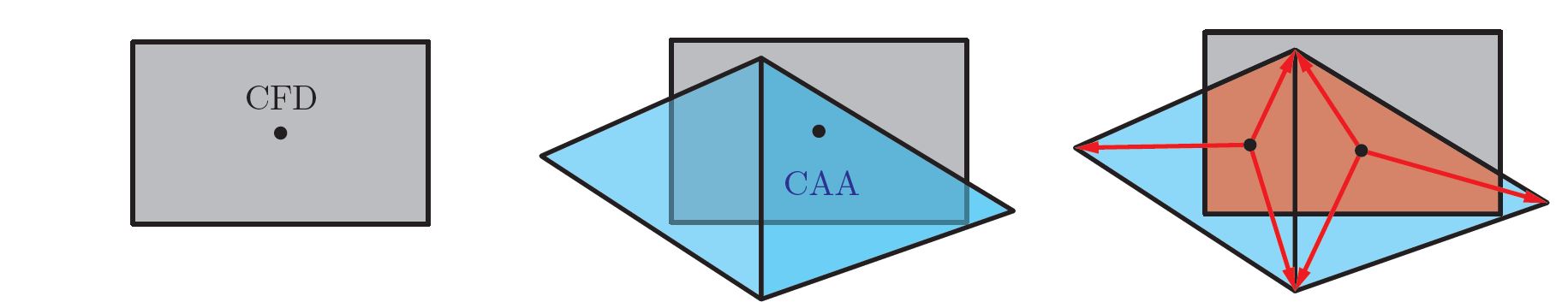}
    \caption{Cut volume cell interpolator.}
    \label{fig:CV}
\end{figure}

The cut volume cell interpolator is defined by the type variable
\begin{lstlisting}[language=XML]
type="FieldInterpolation_Conservative_CutCell"
\end{lstlisting} 
based on figure \ref{fig:CV}, the interpolation is carried out.

\section{Field derivatives}

\subsection{Gradient}

The gradient operator calculates the gradient $\mathbf{r}$ of a scalar source quantity $s$.

$\mathbf{r} = \nabla s$

The derivative is performed on radial basis functions \cite{schoder2020radial}.
The filter is defined as following in the CFSdat xml input:

\begin{lstlisting}[language=XML]
<differentiation type="SpaceDifferentiation_Gradient" inputFilterIds="input" id="gradient">
  <RBF_Settings epsilonScaling="..." betaScaling="..." kScaling="..." logEps="..."/>
  <targetMesh>
	<hdf5 fileName="targetMesh.h5"/>
  </targetMesh>
  <singleResult>
	<inputQuantity resultName="..."/>
	<outputQuantity resultName="..."/>
  </singleResult>
  <regions>
	<sourceRegions>
	  <region name="..." />
	</sourceRegions>
	<targetRegions>
	  <region name="..."/>
	</targetRegions>
  </regions>
</differentiation>
\end{lstlisting} 

In *RBF\_Settings*, the following (optional and mandatory) attributes may be adjusted:
\begin{itemize}
    \item epsilonScaling (mandatory): Controls the "smoothness" of the basis function. The smoother the Gauss-like surface is, the better the results will be BUT only until a certain number, when the matrix becomes so ill-conditioned, which will either result in an exception or very bad results. Typical values ~0.1
    \item betaScaling: Slope of the linear term that was added to the radial basic functions
    \item kScaling: constant term that was added to the radial basic functions
    \item logEps: Console output of [minimal distance, maximal distance, optimized epsilon]. If you run into convergence troubles of the matrix inversion one can activate this. It is only recommended for investigating the quality of the derivative, because it produces a lot of console output.
\end{itemize}

\subsection{Divergence}

The divergence operator calculates the divergence $r$ of a source vector $\mathbf{s}$.
$\mathbf{r} = \nabla \cdot \mathbf{s}$
The filter is defined as following in the CFSdat xml input:

\begin{lstlisting}[language=XML]
<differentiation type="SpaceDifferentiation_Divergence" inputFilterIds="input" id="divergence">
  <RBF_Settings epsilonScaling="..." betaScaling="..." kScaling="..." logEps="..."/>
  <targetMesh>
	<hdf5 fileName="targetMesh.h5"/>
  </targetMesh>
  <singleResult>
	<inputQuantity resultName="..."/>
	<outputQuantity resultName="..."/> 
  </singleResult>
  <regions>
	<sourceRegions>
	  <region name="..." />
	</sourceRegions>
	<targetRegions>
	  <region name="..."/>
	</targetRegions>
  </regions>
</differentiation>
\end{lstlisting}

\subsection{Curl}

The curl operator calculates the curl $\mathbf{r}$ of a source vector $\mathbf{s}$.
$\mathbf{r} = \nabla \times \mathbf{s}$
The filter is defined as following in the CFSdat xml input:
\begin{lstlisting}[language=XML]
<differentiation type="SpaceDifferentiation_Curl" inputFilterIds="input" id="curl">
  <RBF_Settings epsilonScaling="..." betaScaling="..." kScaling="..." logEps="..."/>
  <targetMesh>
	<hdf5 fileName="targetMesh.h5"/>
  </targetMesh>
  <singleResult>
	<inputQuantity resultName="..."/>
	<outputQuantity resultName="..."/>
  </singleResult>
  <regions>
	<sourceRegions>
	  <region name="..."/>
	</sourceRegions>
	<targetRegions>
	  <region name="..."/>
	</targetRegions>
  </regions>
</differentiation>
\end{lstlisting}

\section{Aeroacoustic source terms} 

\subsection{Lamb Vector}
The Lamb Vector filter computes the Lamb vector $\bm L$ based on the velocity $\mathbf{u}$, vorticity $\bm{\omega}$ and the density $\rho$. 

\begin{equation}
\bm L = \bm \omega \times \boldsymbol{u}  ,
\end{equation}

with vorticity as

\begin{equation}
\bm \omega = \nabla \times \boldsymbol{u}  .
\end{equation}

However, it is possible to compute the Lamb vector only based on the velocity, or on the velocity and the vorticity. If only the velocity is defined, the vorticity is computed internally. 

\begin{itemize}
    \item The epsilonScaling Parameter scales the radial basis functions used for computing spatial derivatives \cite{schoder2020radial}. It controls the “smoothness” of the basis function. The smoother the gauss-like surface is, the better the results will be, BUT only until a certain number, when the matrix becomes too ill-conditioned, which will result in very bad results. Typical values: 1e-1 - 1e-4.
    \item The kScaling parameter is an optional parameter and defines a constant term that is added to the radial basic function.
    \item The betaScaling parameter defines the slope of a linear term that is added to the radial basis function.
    \item The logEps Parameter enables a detailes console output (minimal distance, maximal distance, optimized parameters). Therefore, it should only be used if the spatial derivatives are investigated, because it totally spams the console
\end{itemize}

\begin{lstlisting}[language=XML]
    <aeroacoustic type="AeroacousticSource_LambVector" inputFilterIds="..." id="...">
      <RBF_Settings epsilonScaling="1e-4"  kScaling="" betaScaling="" logEps=false/>
      <targetMesh>
        <hdf5 fileName=..."/>
      </targetMesh>
      <ResultList>
        <velocity resultName="..."/>
        <vorticity/>
        <density /> 
        <outputQuantity resultName="..."/>
      </ResultList>
      <regions>
        <sourceRegions>
          <region name="..."/>
        </sourceRegions>
        <targetRegions>
          <region name="..."/>
        </targetRegions>
      </regions>
    </aeroacoustic>
\end{lstlisting} 


\subsection{Lighthill Source Term}

The Lighthill source term filter computes aeroacoustic source terms based on Lighthill's analogy for incompressible flows. Therefore, first the Lighthill source term vector is computed

\begin{equation}
\nabla \cdot \bm T = \nabla(\frac{1}{2} \bm u \cdot \bm u ) + \bm L \, ,
\end{equation}

with $\bm u$ as velocity, and $\bm L$ as Lamb vector, and $\bm T$ as the Lighthill stress tensor. Finally, as actual source term, the divergence of the source term vector is computed and established as outputQuantity. 

\begin{equation}
\nabla \cdot \nabla \cdot \bm T = \nabla \cdot (\nabla(\frac{1}{2} \bm u \cdot \bm u ) + \bm L)
\end{equation}

It is possible to compute the Lighthill source term only based on the velocity, or on the velocity and the vorticity. If only the velocity is defined, the vorticity is computed internally and the parameters are the same as for the Lamb vector.

\begin{lstlisting}[language=XML]
    <aeroacoustic type="AeroacousticSource_LighthillSourceTerm" inputFilterIds="..." id="...">
      <RBF_Settings  epsilonScaling="1e-4"  kScaling="" betaScaling="" logEps=false/>
      <sourceSum>true</sourceSum>
      <targetMesh>
        <hdf5 fileName="..."/>
      </targetMesh>
      <ResultList>
        <velocity resultName="..."/>
        <vorticity/>
        <density/>
        <outputQuantity resultName="..."/>
      </ResultList>
      <regions>
        <sourceRegions>
          <region name="..."/>
        </sourceRegions>
        <targetRegions>
          <region name="..."/>
        </targetRegions>
      </regions>
    </aeroacoustic>
\end{lstlisting}

\subsection{Lighthill Source Term Vector}

The Lighthill source term vector computes a vector corresponding with the $\nabla \cdot \bm T$, where $\bm T$ denotes the Lighthill stress tensor, for incompressible flows.

It is possible to compute the Lighthill source term vector only based on the velocity, or on the velocity and the vorticity. If only the velocity is defined, the vorticity is computed internally and the parameters are the same as for the Lamb vector.

\begin{lstlisting}[language=XML]
    <aeroacoustic type="AeroacousticSource_LighthillSourceTermVector" inputFilterIds="..." id="...">
      <RBF_Settings  epsilonScaling="1e-4"  kScaling="" betaScaling="" logEps=false/>
      <sourceSum>true</sourceSum>
      <targetMesh>
        <hdf5 fileName="..."/>
      </targetMesh>
      <ResultList>
        <velocity resultName="..."/>
        <vorticity/>
        <density/>
        <outputQuantity resultName="..."/>
      </ResultList>
      <regions>
        <sourceRegions>
          <region name="..."/>
        </sourceRegions>
        <targetRegions>
          <region name="..."/>
        </targetRegions>
      </regions>
    </aeroacoustic>
\end{lstlisting}

\subsection{Time Derivative (Simplified PCWE Source Term)}

For low flow velocities, convective effects may be neglected for the PCWE, and the resulting source term simplifies to just the time derivative of the incompressible pressure $\partial p^{\mathrm{ic}}/\partial t$. This is the reason the time derivative filter is placed in the section *Aeroacoustic Source Terms*. Of course, the filter can be applied to a different quantity, of which a time derivative is required, as well.

The time derivative of the desired quantity $\dot{q}(t)$ is computed by a smooth noise-robust differentiator, which suppresses high frequencies and is precise on low frequencies, according to \href{http://www.holoborodko.com/pavel/numerical-methods/numerical-derivative/smooth-low-noise-differentiators/}{this website}. For an efficient and robust calculation the order of the differentiator is set to $N=5$ and the time derivative is calculated by 
\begin{equation}
\dot{q}(t)=\frac{2(q_1-q_{-1})+q_2-q_{-2}}{8 \Delta t},
\end{equation}
where the index of $q$ defines the time step relative to the time step of which the derivative is calculated and $\Delta t$ is the time step size. The computation only requires the definition of the input quantity (*inputQuantity*-tag) and the desired name of the output quantity (*outputQuantity*-tag) as indicated below.

\begin{lstlisting}[language=XML]
<timeDeriv1 id="TimeDerivative" inputFilterIds="input">
	<singleResult>
		<inputQuantity resultName="fluidMechPressure"/>
		<outputQuantity resultName="acouRhsLoadP"/>
	</singleResult>
</timeDeriv1>
\end{lstlisting}

\section{Applications}

In a number of different publications, the details on the module \cite{schoder2018aeroacoustic,schoder2019hybrid,schoder2020conservative} are presented, and applied it to special aeroacoustic source models \cite{kaltenbacher10:2,schoder2020aeroacoustic,schoder2020radial,schoder2019hybrid,weitz2019numerical,schoder2019conservative,schoder2022aeroacoustic,kaltenbacher2021physical}, and special wave equation models \cite{kaltenbacher2018nonconforming,schoder2022cpcwe} are available. Since 2016, the software provided solutions for challenges in acoustical engineering and medicine. Car frame noise \cite{engelmann2020generic,freidhager2021simulationen,schoder2020numerical,weitz2019numerical,maurerlehner2022aeroacoustic}, fan noise \cite{schoder2020computational,tautz2018source,kaltenbacher2017computational,kaltenbacher2016modeling,schoder2021application,tieghi2022machine,tieghi2023machine,schoder2022dataset,schoder2022quantification} noise emissions of the turbocharger compressor \cite{freidhager2022applicability,kaltenbacher2020modelling,freidhager2020influences}, HVAC systems were computed. For post-processing, the fluid field was decomposed into a longitudinal and transversal processes \cite{schoder2020postprocessing,schoder2020postprocessing2,schoder2019helmholtz,schoder2022post}. Furthermore, the human phonation process is studied in detail \cite{schoder2020hybrid,valavsek2019application,zorner2016flow,schoder2021aeroacoustic,falk20213d,lasota2021impact,maurerlehner2021efficient,schoder2022learning,lasota2023anisotropic,schoder2022error,schoder2022pcwe,kraxberger2022machine}. The conservative source term interpolation is used in many aeroacoustic workflows including the FEM-assembly of aeroacoustic source terms of the aeroacoustic wave equation based on Pierce operator \cite{schoder2024aeroacoustic,schoder2023acoustic,schoder2022aeroacoustic} and the compressible perturbed convective wave equation \cite{schoder2022cpcwe,schoder2025perturbed}.

\section{Acknowledge}
We would like to acknowledge the authors of openCFS.

\bibliographystyle{abbrv}
\bibliography{references}  

@article{schoder2020postprocessing,
  title={Postprocessing of Direct Aeroacoustic Simulations Using Helmholtz Decomposition},
  author={Schoder, Stefan and Roppert, Klaus and Kaltenbacher, Manfred},
  journal={AIAA Journal},
  pages={1--9},
  year={2020},
  publisher={American Institute of Aeronautics and Astronautics}
}

@inproceedings{schoder2022aeroacoustic,
  title={Aeroacoustic wave equation based on Pierce's operator applied to the sound generated by a mixing layer},
  author={Schoder, Stefan and Kaltenbacher, Manfred and Spieser, {\'E}tienne and Vincent, Hugo and Bogey, Christophe and Bailly, Christophe},
  booktitle={28th AIAA/CEAS Aeroacoustics 2022 Conference},
  pages={2896},
  year={2022}
}

@article{schoder2023acoustic,
  title={Acoustic modeling using the aeroacoustic wave equation based on Pierce’s operator},
  author={Schoder, Stefan and Spieser, {\'E}tienne and Vincent, Hugo and Bogey, Christophe and Bailly, Christophe},
  journal={AIAA Journal},
  volume={61},
  number={9},
  pages={4008--4017},
  year={2023},
  publisher={American Institute of Aeronautics and Astronautics}
}

@article{schoder2025perturbed,
  title={Perturbed convective wave equation for low-to-medium Mach number subsonic flows},
  author={Schoder, Stefan and Bagheri, Eman and Bogey, Christophe and Bailly, Christophe},
  journal={Journal of Sound and Vibration},
  pages={119549},
  year={2025},
  publisher={Elsevier}
}

@article{schoder2024aeroacoustic,
  title={Aeroacoustic Source Potential Based on Poisson’s Equation},
  author={Schoder, Stefan and Bagheri, Eman and Spieser, {\'E}tienne},
  journal={AIAA Journal},
  volume={62},
  number={7},
  pages={2772--2782},
  year={2024},
  publisher={American Institute of Aeronautics and Astronautics}
}

@article{wurzinger2024pycfs,
  title={pyCFS-data: Data Processing Framework in Python for openCFS},
  author={Wurzinger, Andreas and Schoder, Stefan},
  journal={arXiv preprint arXiv:2405.03437},
  year={2024}
}

@article{tieghi2023machine,
  title={Machine-learning clustering methods applied to detection of noise sources in low-speed axial fan},
  author={Tieghi, Lorenzo and Becker, Stefan and Corsini, Alessandro and Delibra, Giovanni and Schoder, Stefan and Czwielong, Felix},
  journal={Journal of Engineering for Gas Turbines and Power},
  volume={145},
  number={3},
  pages={031020},
  year={2023},
  publisher={American Society of Mechanical Engineers}
}

@article{lasota2023anisotropic,
  title={Anisotropic minimum dissipation subgrid-scale model in hybrid aeroacoustic simulations of human phonation},
  author={Lasota, Martin and {\v{S}}idlof, Petr and Maurerlehner, Paul and Kaltenbacher, Manfred and Schoder, Stefan},
  journal={arXiv preprint arXiv:2301.00606},
  year={2023}
}

@article{CFS,
  title={openCFS: Open Source Finite Element Software for Coupled Field Simulation--Part Acoustics},
  author={Schoder, Stefan and Roppert, Klaus},
  journal={arXiv preprint arXiv:2207.04443},
  year={2022}
}

@inproceedings{freidhager2020influences,
  title={The influences of spatial and temporal discretization in flow simulation on Lighthill’s aeroacoustic source terms applied to a turbocharger},
  author={Freidhager, Clemens and Schoder, Stefan and Kaltenbacher, Manfred},
  booktitle={AIAA AVIATION 2020 FORUM},
  pages={2546},
  year={2020}
}

@article{schoder2020postprocessing2,
  title={Helmholtz’s decomposition for compressible flows and its application to computational aeroacoustics},
  author={Schoder, Stefan and Roppert, Klaus and Kaltenbacher, Manfred},
  journal={SN Partial Differ. Equ. Appl.},
  pages={1--20},
  year={2020},
  publisher={Springer Nature}
}

@ARTICLE{kaltenbacher10:2,
  author = {M.~Kaltenbacher and M.~Escobar and I.~Ali and S.~Becker},
  title = {{Numerical Simulation of Flow-Induced Noise Using LES/SAS and Lighthill's
	Acoustics Analogy}},
  journal = {International Journal for Numerical Methods in Fluids},
  year = {2010},
  volume = {63},
  pages = {1103--1122},
  number = {9}
}

@article{schoder2020conservative,
  title={Conservative interpolation of aeroacoustic sources in a hybrid workflow applied to fan},
  author={Schoder, Stefan J and Junger, Clemens and Weitz, Michael and Kaltenbacher, Manfred},
  journal={arXiv preprint arXiv:2009.02341},
  year={2020}
}

@article{schoder2022dataset,
  title={Dataset FAN-01: Revisiting the EAA Benchmark for a low-pressure axial fan},
  author={Schoder, Stefan and Czwielong, Felix},
  journal={arXiv preprint arXiv:2211.12014},
  year={2022}
}

@article{schoder2022pcwe,
  title={PCWE for FSAI--Derivation of scalar wave equations for fluid-structure-acoustics interaction of low Mach number flows},
  author={Schoder, Stefan},
  journal={arXiv preprint arXiv:2211.07490},
  year={2022}
}

@incollection{kaltenbacher2021physical,
  title={Physical Models for Flow: Acoustic Interaction},
  author={Kaltenbacher, Manfred and Schoder, Stefan},
  booktitle={Waves in Flows},
  pages={265--353},
  year={2021},
  publisher={Springer}
}

@inproceedings{schoder2022quantification,
  title={Quantification of the Acoustic Emissions of an Electric Ducted Fan Unit},
  author={Schoder, Stefan and Schmidt, Jakob and F{\"u}rlinger, Andreas and Kaltenbacher, Manfred},
  booktitle={2022 Delft International Conference on Urban Air-Mobility: DICUAM 2022},
  year={2022}
}

@article{freidhager2022applicability,
  title={Applicability of two hybrid sound prediction methods for assessing in-duct sound absorbers of turbocharger compressors},
  author={Freidhager, Clemens and Schoder, Stefan and Maurerlehner, Paul and Renz, Andreas and Becker, Stefan and Kaltenbacher, Manfred},
  journal={Acta Acustica},
  volume={6},
  pages={37},
  year={2022},
  publisher={EDP Sciences}
}

@inproceedings{schoder2022post,
  title={Post-processing subsonic flows using physics-informed neural networks},
  author={Schoder, Stefan and Museljic, Eniz and Kraxberger, Florian and Wurzinger, Andreas},
  booktitle={2023 AIAA AVIATION Forum},
  year={2022}
}

@article{schoder2022cpcwe,
  title={cPCWE--Perturbed Convective Wave Equation based on Compressible Flows},
  author={Schoder, Stefan},
  journal={arXiv preprint arXiv:2209.11410},
  year={2022}
}

@article{kraxberger2022machine,
  title={Machine-learning applied to classify flow-induced sound parameters from simulated human voice},
  author={Kraxberger, Florian and Wurzinger, Andreas and Schoder, Stefan},
  journal={arXiv preprint arXiv:2207.09265},
  year={2022}
}

@article{maurerlehner2022aeroacoustic,
  title={Aeroacoustic formulations for confined flows based on incompressible flow data},
  author={Maurerlehner, Paul and Schoder, Stefan and Tieber, Johannes and Freidhager, Clemens and Steiner, Helfried and Brenn, G{\"u}nter and Sch{\"a}fer, Karl-Heinz and Ennemoser, Andreas and Kaltenbacher, Manfred},
  journal={Acta Acustica},
  volume={6},
  pages={45},
  year={2022},
  publisher={EDP Sciences}
}

@article{schoder2022error,
  title={Error detection and filtering of incompressible flow simulations for aeroacoustic predictions of human voice},
  author={Schoder, Stefan and Kraxberger, Florian and Falk, Sebastian and Wurzinger, Andreas and Roppert, Klaus and Kniesburges, Stefan and D{\"o}llinger, Michael and Kaltenbacher, Manfred},
  journal={The Journal of the Acoustical Society of America},
  volume={152},
  number={3},
  pages={1425--1436},
  year={2022},
  publisher={Acoustical Society of America}
}

@inproceedings{schoder2019helmholtz,
  title={Helmholtz's Decomposition applied to Aeroacoustics},
  author={Schoder, Stefan and Kaltenbacher, Manfred and Roppert, Klaus},
  booktitle={25th AIAA/CEAS Aeroacoustics Conference},
  year={2019-2561}
}

@article{kaltenbacher2016modeling,
  title={{Modeling and Finite Element Formulation for Acoustic Problems Including Rotating Domains}},
  author={Kaltenbacher, M. and H{\"u}ppe, A. and Grabinger, J. and Wohlmuth, B.},
  journal={AIAA Journal},
  year={2016},
  publisher={American Institute of Aeronautics and Astronautics}
}

@book{cgal,
  title = "{CGAL} User and Reference Manual",
  author = "{The~CGAL~Project}",
  publisher = "{CGAL Editorial Board}",
  edition = "{4.10}",
  year = 2017,
  url = "http://doc.cgal.org/4.10/Manual/packages.html"
}

@book{kaltenbacher2017computational,
  title={{Computational Acoustics}},
  author={Kaltenbacher, M.},
  isbn={9783319590387},
  series={CISM International Centre for Mechanical Sciences},
  year={2017},
  publisher={Springer International Publishing}
}

@techreport{kaltenbacher2020modelling,
  title={Modelling and Numerical Simulation of the Noise Generated by Automotive Turbocharger Compressor},
  author={Kaltenbacher, Manfred and Freidhager, Clemens and Schoder, Stefan},
  year={2020},
  institution={SAE Technical Paper}
}

@inproceedings{schoder2020radial,
  title={Radial Basis Function Interpolation for Computational Aeroacoustics},
  author={Schoder, Stefan and Junger, Clemens and Roppert, Klaus and Kaltenbacher, Manfred},
  booktitle={AIAA AVIATION 2020 FORUM},
  pages={2511},
  year={2020}
}

@article{schoder2019hybrid,
  title={Hybrid Aeroacoustic Computations: State of Art and New Achievements},
  author={Schoder, Stefan and Kaltenbacher, Manfred},
  journal={Journal of Theoretical and Computational Acoustics},
  volume={27},
  number={04},
  pages={1950020},
  year={2019},
  publisher={World Scientific}
}

@article{zorner2016flow,
  title={Flow and acoustic effects in the larynx for varying geometries},
  author={Z{\"o}rner, S and {\v{S}}idlof, P and H{\"u}ppe, A and Kaltenbacher, M},
  journal={Acta Acustica united with Acustica},
  volume={102},
  number={2},
  pages={257--267},
  year={2016},
  publisher={S. Hirzel Verlag}
}

@article{tautz2018source,
  title={Source formulations and boundary treatments for Lighthill’s analogy applied to incompressible flows},
  author={Tautz, Matthias and Besserer, Kerstin and Becker, Stefan and Kaltenbacher, Manfred},
  journal={AIAA Journal},
  volume={56},
  number={7},
  pages={2769--2781},
  year={2018},
  publisher={American Institute of Aeronautics and Astronautics}
}

@article{kaltenbacher2018nonconforming,
  title={Nonconforming Finite Elements Based on Nitsche-Type Mortaring for Inhomogeneous Wave Equation},
  author={Kaltenbacher, Manfred and Floss, Sebastian},
  journal={Journal of theoretical and computational acoustics},
  volume={26},
  number={03},
  pages={1850028},
  year={2018},
  publisher={World Scientific}
}

@article{valavsek2019application,
  title={On the application of acoustic analogies in the numerical simulation of human phonation process},
  author={Val{\'a}{\v{s}}ek, J and Kaltenbacher, M and Sv{\'a}{\v{c}}ek, P},
  journal={Flow, Turbulence and Combustion},
  volume={102},
  number={1},
  pages={129--143},
  year={2019},
  publisher={Springer}
}

@article{weitz2019numerical,
  title={Numerical investigation of the resonance behavior of flow-excited Helmholtz resonators},
  author={Weitz, Michael and Schoder, Stefan and Kaltenbacher, Manfred},
  journal={PAMM},
  volume={19},
  number={1},
  pages={e201900033},
  year={2019},
  publisher={Wiley Online Library}
}

@article{schoder2020computational,
  title={Computational aeroacoustics of the EAA benchmark case of an axial fan},
  author={Schoder, Stefan and Junger, Clemens and Kaltenbacher, Manfred},
  journal={Acta Acustica},
  volume={4},
  number={5},
  pages={22},
  year={2020},
  publisher={EDP Sciences}
}

@article{schoder2020hybrid,
  title={Hybrid aeroacoustic approach for the efficient numerical simulation of human phonation},
  author={Schoder, Stefan and Weitz, Michael and Maurerlehner, Paul and Hauser, Alexander and Falk, Sebastian and Kniesburges, Stefan and D{\"o}llinger, Michael and Kaltenbacher, Manfred},
  journal={The Journal of the Acoustical Society of America},
  volume={147},
  number={2},
  pages={1179--1194},
  year={2020},
  publisher={Acoustical Society of America}
}

@article{schoder2020aeroacoustic,
  title={Aeroacoustic source term computation based on radial basis functions},
  author={Schoder, Stefan and Roppert, Klaus and Weitz, Michael and Junger, Clemens and Kaltenbacher, Manfred},
  journal={International Journal for Numerical Methods in Engineering},
  volume={121},
  number={9},
  pages={2051--2067},
  year={2020},
  publisher={Wiley Online Library}
}

@article{schoder2020numerical,
  title={Numerical investigation of a deep cavity with an overhanging lip considering aeroacoustic feedback mechanism},
  author={Schoder, Stefan and Lazarov, Ivan and Kaltenbacher, Manfred},
  journal={arXiv preprint arXiv:2006.03279},
  year={2020}
}

@techreport{engelmann2020generic,
  title={A Generic Testbody for Low-Frequency Aeroacoustic Buffeting},
  author={Engelmann, Rafael and Gabriel, Christoph and Schoder, Stefan and Kaltenbacher, Manfred},
  year={2020},
  institution={SAE Technical Paper}
}

@phdthesis{schoder2018aeroacoustic,
  title={Aeroacoustic analogies based on compressible flow data},
  author={Schoder, Stefan Josef},
  year={2018},
  school={Wien}
}

@article{schoder2022learning,
  title={Learning Expertise Actively to model Domain knowledge (LEAD) with application to human phonation},
  author={Schoder, Stefan and Roppert, Klaus},
  year={2022},
  journal={arXiv}
}

@inproceedings{tieghi2022machine,
  title={Machine-Learning Clustering Methods Applied to Detection of Noise Sources in Low-Speed Axial Fan},
  author={Tieghi, Lorenzo and Becker, Stefan and Corsini, Alessandro and Delibra, Giovanni and Schoder, Stefan and Czwielong, Felix},
  booktitle={2022 Turbomachinery Technical Conference \& Exposition: ASME Turbo Expo 2022},
  year={2022}
}

@article{maurerlehner2021efficient,
  title={Efficient numerical simulation of the human voice},
  author={Maurerlehner, Paul and Schoder, Stefan and Freidhager, Clemens and Wurzinger, Andreas and Hauser, Alexander and Kraxberger, Florian and Falk, Sebastian and Kniesburges, Stefan and Echternach, Matthias and D{\"o}llinger, Michael and others},
  journal={e \& i Elektrotechnik und Informationstechnik},
  volume={138},
  number={3},
  pages={219--228},
  year={2021},
  publisher={Springer}
}

@article{lasota2021impact,
  title={Impact of the sub-grid scale turbulence model in aeroacoustic simulation of human voice},
  author={Lasota, Martin and {\v{S}}idlof, Petr and Kaltenbacher, Manfred and Schoder, Stefan},
  journal={Applied Sciences},
  volume={11},
  number={4},
  pages={1970},
  year={2021},
  publisher={MDPI}
}

@inproceedings{schoder2019conservative,
  title={Conservative source term interpolation for hybrid aeroacoustic computations},
  author={Schoder, Stefan and Junger, Clemens and Weitz, Michael and Kaltenbacher, Manfred},
  booktitle={25th AIAA/CEAS aeroacoustics conference},
  pages={2538},
  year={2019}
}

@article{freidhager2021simulationen,
  title={Simulationen von Str{\"o}mungsakustik in rotierenden Bauteilen zur Entwicklung von Antriebskonzepten der Autos der Zukunft},
  author={Freidhager, Clemens and Maurerlehner, Paul and Roppert, Klaus and Wurzinger, Andreas and Hauser, Alexander and Heinisch, Martin and Schoder, Stefan and Kaltenbacher, Manfred},
  journal={e \& i Elektrotechnik und Informationstechnik},
  volume={138},
  number={3},
  pages={212--218},
  year={2021},
  publisher={Springer}
}

@article{falk20213d,
  title={3D-FV-FE aeroacoustic larynx model for investigation of functional based voice disorders},
  author={Falk, Sebastian and Kniesburges, Stefan and Schoder, Stefan and Jakuba{\ss}, Bernhard and Maurerlehner, Paul and Echternach, Matthias and Kaltenbacher, Manfred and D{\"o}llinger, Michael},
  journal={Frontiers in physiology},
  volume={12},
  pages={616985},
  year={2021},
  publisher={Frontiers Media SA}
}

@article{schoder2021aeroacoustic,
  title={Aeroacoustic sound source characterization of the human voice production-perturbed convective wave equation},
  author={Schoder, Stefan and Maurerlehner, Paul and Wurzinger, Andreas and Hauser, Alexander and Falk, Sebastian and Kniesburges, Stefan and D{\"o}llinger, Michael and Kaltenbacher, Manfred},
  journal={Applied Sciences},
  volume={11},
  number={6},
  pages={2614},
  year={2021},
  publisher={MDPI}
}

@article{schoder2021application,
  title={Application limits of conservative source interpolation methods using a low Mach number hybrid aeroacoustic workflow},
  author={Schoder, Stefan and Wurzinger, Andreas and Junger, Clemens and Weitz, Michael and Freidhager, Clemens and Roppert, Klaus and Kaltenbacher, Manfred},
  journal={Journal of Theoretical and Computational Acoustics},
  volume={29},
  number={01},
  pages={2050032},
  year={2021},
  publisher={World Scientific}
}






\end{document}